\documentclass{amsart}
\usepackage{amssymb,amsmath,latexsym,amsthm}
\usepackage[english]{babel}

\usepackage{amsmath,amssymb,amsbsy,amsfonts,amsthm,latexsym,
                       amsopn,amstext,amsxtra,euscript,amscd}

\usepackage{array}
\usepackage{ifthen}
\usepackage{url}
\usepackage{hyperref}
            
\newtheorem{lem}{Lemma}
\newtheorem{lemma}[lem]{Lemma}

\newtheorem{thm}{Theorem}
\newtheorem{theorem}[thm]{Theorem}

\def\\{\cr}
\def\({\left(}
\def\){\right)}
\def\[{\left[}
\def\]{\right]}
\def\<{\langle}
\def\>{\rangle}

\begin{document}

\title[On a nonintegrality conjecture]{On a nonintegrality conjecture}

\author{Florian Luca}
\address{School of Maths, Wits University, South Africa, Max Planck Institute for Mathematics, Bonn, Germany, Research Group in Algebraic Structures and Applications, King Abdulaziz University, Jeddah, Saudi Arabia}
\email{florian.luca@wits.ac.za}

\author{Carl Pomerance}
\address{Mathematics Department, Dartmouth College, 27 N. Main Street, 6188 Kemeny Hall, Hanover, Hanover NH 03755-3551, USA}
\email{carlp@math.dartmouth.edu}

\date{\today}

\pagenumbering{arabic}

\begin{abstract}
It is conjectured that the sum
$$
S_r(n)=\sum_{k=1}^{n} \frac{k}{k+r}\binom{n}{k}
$$
for positive integers $r,n$ is never integral.  This has been shown for
$r\le 22$.  In this note we study the problem in the ``$n$ aspect"
showing that the set of $n$ such that $S_r(n)\in {\mathbb Z}$ for some $r\ge 1$ has asymptotic density $0$.  Our principal tools are some deep results on the
distribution of primes in short intervals.
\end{abstract}



\maketitle

\section{Introduction}

For positive integers $r,n$ let
$$
S_r(n)=\sum_{k=1}^n \frac{k}{k+r}\binom{n}{k}.
$$
Motivated by some cases with small $r$, L\'opez-Aguayo \cite{L}
asked if $S_r(n)$ is ever an integer, showing for $r\in \{1,2,3,4\}$
that $S_r(n)$ is not integral for all $n$.   In \cite{LL} it was conjectured that $S_r(n)$ is never
integral, and they proved the conjecture for $r\le 6$.
In \cite{LLPT} it was proved for $r\le 22$.  Also in \cite{LLPT}, using a deep
theorem of Montgomery and Vaughan \cite{MV}, it was shown  for a fixed $r$ 
that the set of $n$ such that $S_r(n)\in {\mathbb Z}$ has upper density bounded by $O_k(1/r^k)$ 
for any $k\ge 1$.  In fact, this density is 0, as we shall show.  Actually we prove
a stronger result.
Let 
$$
{\mathcal S}:=\{n: S_r(n)\in {\mathbb Z}~{\text{\rm for~some}}~r\ge 1\}.
$$

\begin{theorem}
The set ${\mathcal S}$ has zero density as a subset of the integers. 
\end{theorem}

It follows from our argument that if we put ${\mathcal S}(x)={\mathcal S}\cap [1,x]$ then $\#{\mathcal S}(x)=O_A( x/(\log x)^A)$ for every fixed $A$.  In particular,
taking $A=2$, we see that the reciprocal sum of $\mathcal S$ is finite.

\section{The proof}

We let $x$ be large and $n\in {\mathcal S}\cap [x/2,x)$. Thus, 
$S_r(n)\in {\mathbb Z}$ for some $r\ge 1$.  Let
$$
S(r,n):=\sum_{k=0}^n \frac{r}{k+r}\binom{n}{k},
$$
so that $S(r,n)+S_r(n)=\sum_{k=0}^n\binom nk=2^n\in\mathbb Z$, so that
$S(r,n)\in\mathbb Z$.  It is shown in \cite{LL} that
\begin{equation}
\label{eq:2power}
S(r,n)=\sum_{j=1}^r(-1)^{r-j}r\binom{r-1}{j-1}\frac{2^{n+j}-1}{n+j}.
\end{equation}
\begin{lemma}
\label{lem:easy}
If there is a prime $p>n$ that divides one of $1+r,2+r,\dots,n+r$,
then $S_r(n)$ is not integral.
\end{lemma}
\begin{proof}
Write $p$ as $k_0+r$, where $1\le k_0\le n$.  Since $p>n$, we have that
$p$ does not divide any other $k+r$ for $1\le k\le n$.  So the term $(k_0/(k_0+r))\binom n{k_0}$ in the definition
of $S_r(n)$, in reduced form, has a factor $p$ in the denominator, and no other
terms $(k/(k+r))\binom nk$ have this property. 
We deduce that $S_r(n)$ is
nonintegral, completing the proof.
\end{proof}

We distinguish various cases.

\medskip

{\bf Case 1.} $r\ge n$. 

\medskip 

By Sylvester's theorem, one of the integers $1+r,2+r,\dots,n+r$
 is divisible by a prime $p>n$.  It follows 
from Lemma \ref{lem:easy} that $S_r(n)$ is nonintegral.
  From now on, we assume that $n>r$. 

\medskip

{\bf Case 2.} $n>r>(x/2)^{1/10}$. 

\medskip

By a result of Jia (see \cite{Jia}) for every fixed $\varepsilon>0$, the interval 
$[n+1,n+n^{1/20+\varepsilon}]$ contains a prime number $p$ for almost all $n$,
with the number of exceptional values of $n\le x$ being $\ll_{\epsilon,A}x/(\log x)^A$
for every fixed $A>0$.
 If $r>(x/2)^{1/10}\ge(n/2)^{1/10}$, then $r>n^{1/11}$ holds for all $x>x_0$.  
 If   $n$ is not exceptional in the sense of Jia's theorem, then the interval $[n,n+r]$ 
 contains the interval $[n+1,n+n^{1/11}]$ and hence a prime $p>n>r$, so $S_r(n)$ 
 cannot be an integer by Lemma \ref{lem:easy}.  Hence, $n$ must be exceptional in 
 the sense of Jia's theorem and the set of such $n$ has counting function 
$O_A(x/(\log x)^A)$ for any fixed $A>0$. 

\medskip

{\bf Case 3.} $y\le r\le (x/2)^{1/10}$, where $y:=x^{1/\log\log x}$.

\medskip

This is the most interesting part. We prove the following lemma. For an odd prime $p$ we write $\ell_2(p)$ for the order of $2$ modulo $p$.

\begin{lemma}
There exists $r_0$ such that if $r>r_0$, then the interval $I=[r,r+r^{0.61}]$ contains $6$ primes $p_1,\dots,p_6$ such that each $\ell_2(p_i)>r^{0.3}$ for
$1\le i\le6$ and each $\gcd(p_i-1,p_j-1)<r^{0.001}$ for $1\le i<j\le 6$.
\end{lemma}
\begin{proof}
Let $\pi(I)$ be the number of primes in $I$.  From  Baker, Harman, and Pintz
\cite{BHP} we have for large $r$ that
$$
\pi(I)\gg r^{0.61}/\log r.
$$
(Actually, this follows from earlier results, but \cite{BHP} holds the record
currently for primes in short intervals.)
Let ${\mathcal Q}$ be the subset of primes $p\in I$ such that $\ell_2(p)\le r^{0.3}$. By a classical argument, $\#{\mathcal Q}\ll r^{0.6}/\log r$. Indeed,  
$$
r^{\#{\mathcal Q}}\le \prod_{p\in {\mathcal Q}} p\le \prod_{t\le r^{0.3}} (2^t-1)<2^{\sum_{t\le r^{0.3}} t}<2^{r^{0.6}},
$$
from which we deduce the desired upper bound on $\#{\mathcal Q}$. Since 
$$
r^{0.6}/\log r=o(r^{0.61}/\log r)=o(\pi(I)),\qquad {\text{\rm as}}\qquad r\to\infty,
$$ 
we deduce that most primes $p$ in $I$ have $\ell_2(p)\ge r^{0.3}$.  Let
$\mathcal P$ denote this set of primes in $I$, so that $\#{\mathcal P}\gg r^{0.61}/\log r$.
  For any positive integer $d$
the number of pairs of primes $p,q$ in $\mathcal P$ with $d\mid p-1$
and $d\mid q-1$ is $\ll r^{2\times0.61}/d^2$ even ignoring the primality
condition.  Summing over $d\ge r^{0.001}$ we see that the number of pairs
$p,q\in\mathcal P$ with $\gcd(p-1,q-1)\ge r^{0.001}$ is 
$\ll r^{2\times0.61-0.001}$, so that most pairs of
primes $p,q\in \mathcal P$ have $\gcd(p-1,q-1)<r^{0.001}$.  In fact, the
number of 6-tuples of primes $p_1,\dots,p_6\in\mathcal P$ with
some $\gcd(p_i-1,p_j-1)\ge r^{0.001}$ is $\ll r^{6\times0.61-0.001}$,
 so we may deduce that most 6-tuples of primes in $\mathcal P$
satisfy the gcd condition of the lemma.  Of course ``6" may be replaced with
any fixed positive integer, only affecting the choice of $r_0$.
\end{proof}

Let $\{p_1,\ldots,p_6\}$ be the $6$ primes in $I$ which exist for $x>x_0$ (such that $y>r_0$). 
Either there are 4 of these primes such that the interval $[n+1,n+r]$ contains
a multiple of each, or there are 3 of these primes which do not have multiples
in $[n+1,n+r]$.   Take the case of 4 of the primes having a multiple in $[n+1,n+r]$
and without essential loss of generality, say they are $p_1,p_2,p_3,p_4$.
They determine integers $j_1,j_2,j_3,j_4$ with $1\le j_i\le r$ and
$p_i\mid n+j_i$.  However, there is another restriction on $n$ caused
by $S(r,n)$ being integral.   We have each
$\ell_2(p_i)\mid n+j_i$, since otherwise the $j_i$ term in \eqref{eq:2power}
in reduced form contains a factor of $p_i$ in the denominator, a property
not shared with any other term.  This would imply that $S(r,n)$ is nonintegral,
a contradiction.  Thus, we have $\ell_2(p_i)\mid n+j_i$ as claimed for $i=1,2,3,4$.
We conclude that $n$ is in a residue class modulo
$$
M:={\rm lcm}\{p_1,p_2,p_3,p_4,\ell_2(p_1),\ell_2(p_2),\ell_2(p_3),\ell_2(p_4)\}.
$$
Now $p_1,p_2,p_3,p_4$ are distinct primes in $[r+1,r+r^{0.61}]$, and
each $\ell_2(p_i)$, since it divides $p_i-1$, has
all prime factors $\le r$, so is coprime to the other $p_j$'s.  Moreover,
each $\ell_2(p_i)>r^{0.3}$ and being a divisor of $p_i-1$, each
$\gcd(\ell_2(p_i),\ell_2(p_j))\le r^{0.001}$.  Thus,
$$
M>r^4r^{1.2}r^{-0.006}=r^{5.194}.
$$
Further, $M\ll r^8<x$.
Thus, the number of $n$ in this residue class is $\ll x/M<x/r^{5.194}$.  Summing
over the different possibilities for $j_1,j_2,j_3,j_4$, our count is
$\ll x/r^{1.194}$.  Now summing over $r>y$, we have that the number
of $n$ in this case is $\ll x/y^{0.194}$.

We also must consider the possibility that 3 of our 6 primes do not divide
any $n+j$ with $1\le j\le r$.  Again without essential loss of generality,
assume they are $p_1,p_2,p_3$.  Since each is in $[r+1,r+r^{0.61}]$, it
follows that each $p_i$ corners $n$ in a set of $O(r^{0.61})$ residue
classes mod $p_i$.
With the Chinese Remainder Theorem, such $n$'s are in a set of $O(r^{1.83})$ residues classes
modulo $p_{1}p_{2}p_{3}$. Note that the modulus is small, at most $O(r^3)=o(x)$. Thus, the number of such $n$ is at most
$$
O\left(\frac{r^{1.83} x}{p_{1}p_{2}p_{3}}\right)=O\left(\frac{x}{r^{1.17}}\right).
$$
Varying the 3 primes in $\binom{6}{3}=20$ ways multiplies the above count by a constant factor.  Summing on $r>y$ we deduce that the number of $n$
in $(x/2,x]$ is $\ll x/y^{0.17}$.  With our above estimate, this puts the
count in Case 3 at $O(x/y^{0.17})=o(x)$ as $x\to\infty$.

\medskip

{\bf Case 4.} {\it We assume that $r\in (22,y]$.}

\medskip

Here, we do the ``regular" thing, where we distinguish between smooth
numbers and numbers with a large prime factor.  Let $P(m)$ denote the
largest prime factor of $m$.  If $P(n+1)\le y$, this puts $n$ in a set of
size $x/(\log x)^{(1+o(1))\log\log\log x}$ as $x\to\infty$, by standard
estimates for smooth numbers.  So, assume that $p=P(n+1)>y$.
Since $r\le y$, it follows that $p$ does not divide any other $n+j$ with
$j\le r$, so that \eqref{eq:2power} and $S(r,n)$ integral imply that
$\ell_2(p)\mid n+1$.

The number of primes $2<q\le t$ with $\ell_2(q)\le q^{0.3}$ is by
the argument in the previous case  at most $ t^{0.6}$.  By a partial
summation argument, the number of $n\in(x/2,x]$ with $n+1$ divisible by such
a prime $q>y$ is $O(x/y^{0.4})$.  So, assume that $\ell_2(p)>p^{0.3}$.
 The number of integers $n\in(x/2,x]$ with $n+1$ divisible by $p\ell_2(p)$
 is $\le x/(p\ell_2(p))\le x/p^{1.3}$.  Summing on $p>y$ our count is
 $\ll x/y^{0.3}$.
 
\medskip

Putting together everything, we get that  $\#{\mathcal S}(x)$ is
$O_A(x/(\log x)^A)$ for every fixed $A>0$.  This completes the proof of
the theorem.

\medskip
\noindent{\bf Remarks}.
Note that assuming Cram\'er's conjecture that for some constant $c$
and for large $x$
there is a prime in $[x,x+c(\log x)^2]$, the estimate in Case 2 is eliminated.
By then optimizing the choice of $y$, our final count for $\mathcal S(x)$ would
be of the shape $O(x/\exp(c\sqrt{\log x\log\log x}))$ for some $c>0$.  The
hardest cases to try and do better seem to be $r=O(1)$.

Let $s_r(m)$ be the largest $r$-smooth divisor of $m$ and let
$M_r(n)=\min\{s_r(n+j):1\le j\le r\}$.  It follows from \cite[Proposition 3.1]{LLPT}
that if $M_r(n)\le\log_2r$, then $S_r(n)$ is nonintegral.  Unfortunately,
as discussed in \cite[Remark 2]{LLPT}, it is not always the case that
$M_r(n)\le\log_2r$.  Nevertheless, it seems interesting to
get estimates for $M(r):=\max\{M_r(n):n>0\}$.

\end{document}